\documentclass[12pt]{article}
\usepackage{amssymb}
\usepackage{mathrsfs}
\usepackage{pst-poly}     
\usepackage{cite}
\usepackage{amsmath}

\newtheorem{thm}{Theorem}[section]
\newtheorem{lem}[thm]{Lemma}
\newtheorem{Def}[thm]{Definition}
\newtheorem{cor}[thm]{Corollary}

\newenvironment{pf}[1][Proof]{\noindent\textbf{#1.} }{\hfill\rule{1mm}{2mm}}

\makeatletter \@addtoreset{equation}{section} \makeatother

\begin{document}
\title{\bf Generalized Measures of Fault
Tolerance in Exchanged Hypercubes \thanks {The work was supported by
NNSF of China (No.11071233, 61272008).}}
\author
{Xiang-Jun Li$^{a,b}$
\quad Jun-Ming Xu$^a$
\footnote{Corresponding author: xujm@ustc.edu.cn (J.-M. Xu)}\\
{\small $^a$School of Mathematical Sciences}\\
{\small  University of Science and Technology of China}\\
{\small  Wentsun Wu Key Laboratory of CAS, Hefei, 230026, China}  \\
{\small $^b$School of Information and Mathematics}\\
{\small  Yangtze University, Jingzhou, Hubei, 434023, China}\\
 }
\date{}
 \maketitle

\begin{abstract}

The exchanged hypercube $EH(s,t)$, proposed by Loh {\it et al.} [The
exchanged hypercube, IEEE Transactions on Parallel and Distributed
Systems 16 (9) (2005) 866-874], is obtained by removing edges from a
hypercube $Q_{s+t+1}$. This paper considers a kind of generalized
measures $\kappa^{(h)}$ and $\lambda^{(h)}$ of fault tolerance in
$EH(s,t)$ with $1\leqslant s\leqslant t$ and determines
$\kappa^{(h)}(EH(s,t))=\lambda^{(h)}(EH(s,t))= 2^h(s+1-h)$ for any
$h$ with $0\leqslant h\leqslant s$. The results show that at least
$2^h(s+1-h)$ vertices (resp. $2^h(s+1-h)$ edges) of $EH(s,t)$ have
to be removed to get a disconnected graph that contains no vertices
of degree less than $h$, and generalizes some known results.

\vskip6pt


\end{abstract}

\section{Introduction}
It is well known that interconnection networks play an important
role in parallel computing/communication systems. An interconnection
network can be modeled by a graph $G=(V, E)$, where $V$ is the set
of processors and $E$ is the set of communication links in the
network. For graph terminology and notation not defined here we
follow \cite{xu01}.

A subset $S\subset V(G)$ (resp. $F\subset E(G)$) of a connected
graph $G$ is called a {\it vertex-cut} (resp. {\it edge-cut}) if
$G-S$ (resp. $G-F$) is disconnected. The {\it connectivity}
$\kappa(G)$ (resp. {\it edge-connectivity} $\lambda(G)$ ) of $G$ is
defined as the minimum cardinality over all vertex-cuts (resp.
edge-cuts) of $G$. The connectivity $\kappa(G)$ and
edge-connectivity $\lambda(G)$ of a graph $G$ are two important
measurements for fault tolerance of the network since the larger
$\kappa(G)$ or $\lambda(G)$ is, the more reliable the network is.

Because the connectivity has some shortcomings,
Esfahanian~\cite{E89} proposed the concept of restricted
connectivity, Latifi {\it et al.}~\cite{LHP94} generalized it to
restricted $h$-connectivity which can measure fault tolerance of an
interconnection network more accurately than the classical
connectivity. The concepts stated here are slightly different from
theirs.

A subset $S\subset V(G)$ (resp. $F\subset E(G)$) of a connected
graph $G$, if any, is called an {\it $h$-vertex-cut} (resp. {\it
edge-cut}), if $G-S$ (resp. $G-F$) is disconnected and has the
minimum degree at least $h$. The {\it $h$-connectivity} (resp. {\it
edge-connectivity}) of $G$, denoted by $\kappa^{(h)}(G)$ (resp.
$\lambda^{(h)}(G)$), is defined as the minimum cardinality over all
$h$-vertex-cuts (resp. $h$-edge-cut) of $G$. It is clear that, for
$h\geqslant 1$, if $\kappa^{(h)}(G)$ and $\lambda^{(h)}(G)$) exists,
then $\kappa^{(h-1)}(G)\leqslant \kappa^{(h)}(G)$ and
$\lambda^{(h-1)}(G)\leqslant \lambda^{(h)}(G)$. For any graph $G$
and any integer $h$, determining $\kappa^{(h)}(G)$ and
$\lambda^{(h)}(G)$ is quite difficult. In fact, the existence of
$\kappa^{(h)}(G)$ and $\lambda^{(h)}(G)$ is an open problem so far
when $h\geqslant 1$. Only a little knowledge of results have been
known on $\kappa^{(h)}$ and $\lambda^{(h)}$ for particular classes
of graphs and small $h$'s, such as
~\cite{EH88,lx12a,lx12,m10,m11,oc93,wz09,wg98,x00c,xxz05,ylg10,zxh07,zxl06}.

It is widely known that the hypercube $Q_n$ has been one of the most
popular interconnection networks for parallel computer/communication
system. Xu~\cite{x00c} determined $\lambda^{(h)}(Q_n)=2^h(n-h)$ for $h\leqslant n-1$, and
Oh {\it et al.}~\cite{oc93} and Wu {\it et al.}~\cite{wg98}
independently determined $\kappa ^{(h)}(Q_n)=2^h(n-h)$ for
$h\leqslant n-2$.

This paper is concerned about the exchanged hypercubes $EH(s,t)$,
proposed by Loh {\it et al.}~\cite{lhp05}. As a variant of the
hypercube, $EH(s,t)$ is a graph obtained by removing edges from a
hypercube $Q_{s+t+1}$. It not only keeps numerous desirable
properties of the hypercube, but also reduced the interconnection
complexity. Very recently, Ma {\it et al.}~\cite{m11} have
determined $\kappa^{(1)}(EH(s,t))=\lambda^{(1)}(EH(s,t))=2s$. We, in
this paper, will generalize this result by proving that $\kappa
^{(h)}(EH(s,t))=\lambda^{(h)}(EH(s,t))=2^h(s+1-h)$ for any $h$ with
$0\leqslant h\leqslant s$.

The proof of this result is in Section 3. In Section 2, we recall
the structure of $EH(s,t)$ and some lemmas used in our proofs.

\section{Definitions and lemmas}

For a given position integer $n$, let $I_{n}=\{1,2,\ldots, n\}$. The
sequence $x_{n}x_{n-1}\cdots x_{1}$ is said a binary string of
length $n$ if $x_r\in\{0,1\}$ for each $r\in I_n$. Let
$x=x_{n}x_{n-1}\cdots x_{1}$ and $y=y_{n}y_{n-1}\cdots y_{1}$ be two
distinct binary string of length $n$. {\it Hamming distance} between
$x$ and $y$, denoted by $H(x,y)$, is the number of $r$'s for which
$|x_r-y_r|=1$ for $r\in I_n$.

For a binary string $u=u_{n}u_{n-1}\cdots u_{1}u_0$ of length $n+1$,
we call $u_r$ the $r$-th bit of $u$ for $r\in I_n$, and $u_0$ the
last bit of $u$, denote sub-sequence $u_ju_{j-1}\cdots u_{i+1}u_i$
of $u$ by $u[j: i]$, i.e., $u[j,i]=u_ju_{j-1}\cdots u_{i+1}u_i$. Let
$$
 V(s,t) = \{u_{s+t}\cdots
 u_{t+1}u_{t}\cdots u_{1}u_0|\ u_0,u_i\in \{0,1\}, i\in
 I_{s+t}\}.
 $$

\begin{Def}\label{def2.1} The exchanged hypercube is an undirected
graph $EH(s,t)=(V, E)$, where $s\geqslant 1$ and  $t\geqslant 1$ are
integers. The set of vertices $V$ is $V(s,t)$, and the set of edges
$E$ is composed of three disjoint types $E_1, E_2$ and $E_3$.
 $$
 \begin{array}{rl}
 E_1 = &\{u v \in V\times V|\  u[s+t:1] = v[s+t:1],
u_0\not = v_0\},\\
E_2 =&\{ u v \in V\times V|\  u[s +t : t + 1] = v[s +t : t + 1],\\
&   H(u[t : 1], v[t : 1]) = 1,
u_0=v_0=1\},\\
E_3 =& \{u v \in V\times V|\   u[t : 1] = v[t : 1],\\
&  H(u[s +t : t + 1], v[s +t : t +1]) = 1, u_0=v_0=0\}.
\end{array}$$
\end{Def}

Now we give an alternative definition of $EH(s,t)$.

 \begin{Def}\label{def2.2}
An exchanged hypercube $EH(s,t)$ consists of the vertex-set $V(s,t)$
and the edge-set $E$,
two vertex $u = u_{s+t}\cdots u_{t+1}u_{t}\cdots u_{1}u_0$ and $v =
v_{s+t}\cdots v_{t+1}v_{t}\cdots v_{1}v_0$ linked by an edge, called
$r$-dimensional edge, if and only if the
following conditions are satisfied:\\
a). $u$ and $v$ differ exactly in one bit on the $r$-th bit or on the last bit. \\
b). if $r\in I_{t}$,  then $u_0=v_0=1$,\\
c). if $r\in I_{s+t}-I_{t}$,  then $u_0=v_0=0$.
\end{Def}

\begin{figure}[h]
\begin{center}
\begin{pspicture}(0,-2.3429687)(11.162812,2.3429687)
\psdots[dotsize=0.12](5.7409377,1.8045312)
\psdots[dotsize=0.12](7.1809373,1.8045312)
\psdots[dotsize=0.12](5.7409377,0.6045312)
\psdots[dotsize=0.12](7.1809373,0.6045312)
\psdots[dotsize=0.12](5.7409377,-0.59546876)
\psdots[dotsize=0.12](7.1809373,-0.59546876)
\psdots[dotsize=0.12](5.7409377,-1.7954688)
\psdots[dotsize=0.12](7.1809373,-1.7954688)
\psline[linewidth=0.04cm](5.7409377,1.8045312)(5.7409377,0.6045312)
\psline[linewidth=0.04cm](5.7409377,1.8045312)(7.1809373,1.8045312)
\psline[linewidth=0.04cm](7.1809373,1.8045312)(7.1809373,0.6045312)
\psline[linewidth=0.04cm](7.1809373,0.6045312)(5.7409377,-0.59546876)
\psline[linewidth=0.04cm](5.7409377,-0.59546876)(5.7409377,-1.7954688)
\psline[linewidth=0.04cm](5.7409377,-1.7954688)(7.1809373,-1.7954688)
\psline[linewidth=0.04cm](7.1809373,-1.7954688)(7.1809373,-0.59546876)
\psline[linewidth=0.04cm](5.7409377,0.6045312)(7.1809373,-0.59546876)
\usefont{T1}{ptm}{m}{n}
\rput(7.222344,2.1545312){\scriptsize$0001$}
\usefont{T1}{ptm}{m}{n}
\rput(7.3623437,0.23453125){\scriptsize$0101$}
\usefont{T1}{ptm}{m}{n}
\rput(7.3623437,-0.24546875){\scriptsize$1001$}
\usefont{T1}{ptm}{m}{n}
\rput(7.222344,-2.1654687){\scriptsize$1101$}
\usefont{T1}{ptm}{m}{n}
\rput(5.2,1.9145312){\scriptsize$0000$}
\usefont{T1}{ptm}{m}{n}
\rput(5.2,0.47453126){\scriptsize$1000$}
\usefont{T1}{ptm}{m}{n}
\rput(5.2,-0.48546875){\scriptsize$0100$}
\usefont{T1}{ptm}{m}{n}
\rput(5.2,-1.6854688){\scriptsize$1100$}
\psdots[dotsize=0.12](8.620937,1.8045312)
\psdots[dotsize=0.12](10.060938,1.8045312)
\psdots[dotsize=0.12](8.620937,0.6045312)
\psdots[dotsize=0.12](10.060938,0.6045312)
\psdots[dotsize=0.12](8.620937,-0.59546876)
\psdots[dotsize=0.12](10.060938,-0.59546876)
\psdots[dotsize=0.12](8.620937,-1.7954688)
\psdots[dotsize=0.12](10.060938,-1.7954688)
\psline[linewidth=0.04cm](8.620937,1.8045312)(8.620937,0.6045312)
\psline[linewidth=0.04cm](8.620937,1.8045312)(10.060938,1.8045312)
\psline[linewidth=0.04cm](10.060938,1.8045312)(10.060938,0.6045312)
\psline[linewidth=0.04cm](10.060938,0.6045312)(8.620937,-0.59546876)
\psline[linewidth=0.04cm](8.620937,-0.59546876)(8.620937,-1.7954688)
\psline[linewidth=0.04cm](8.620937,-1.7954688)(10.060938,-1.7954688)
\psline[linewidth=0.04cm](10.060938,-1.7954688)(10.060938,-0.59546876)
\psline[linewidth=0.04cm](8.620937,0.6045312)(10.060938,-0.59546876)
\usefont{T1}{ptm}{m}{n}
\rput(10.582344,1.9145312){\scriptsize$0010$}
\usefont{T1}{ptm}{m}{n}
\rput(10.582344,0.47453126){\scriptsize$1010$}
\usefont{T1}{ptm}{m}{n}
\rput(10.582344,-0.48546875){\scriptsize$0110$}
\usefont{T1}{ptm}{m}{n}
\rput(10.582344,-1.6854688){\scriptsize$1110$}
\usefont{T1}{ptm}{m}{n}
\rput(8.6,2.1545312){\scriptsize$0011$}
\usefont{T1}{ptm}{m}{n}
\rput(8.422344,0.23453125){\scriptsize$0111$}
\usefont{T1}{ptm}{m}{n}
\rput(8.422344,-0.24546875){\scriptsize$1011$}
\usefont{T1}{ptm}{m}{n}
\rput(8.6,-2.1654687){\scriptsize$1111$}
\psdots[dotsize=0.12](0.9409375,1.8045312)
\psdots[dotsize=0.12](2.3809376,1.8045312)
\psdots[dotsize=0.12](0.9409375,0.6045312)
\psdots[dotsize=0.12](2.3809376,0.6045312)
\psdots[dotsize=0.12](0.9409375,-0.59546876)
\psdots[dotsize=0.12](2.3809376,-0.59546876)
\psdots[dotsize=0.12](0.9409375,-1.7954688)
\psdots[dotsize=0.12](2.3809376,-1.7954688)
\psline[linewidth=0.04cm](0.9409375,1.8045312)(0.9409375,0.6045312)
\psline[linewidth=0.04cm](0.9409375,1.8045312)(2.3809376,1.8045312)
\psline[linewidth=0.04cm](2.3809376,1.8045312)(2.3809376,0.6045312)
\psline[linewidth=0.04cm](2.3809376,0.6045312)(0.9409375,-0.59546876)
\psline[linewidth=0.04cm](0.9409375,-0.59546876)(0.9409375,-1.7954688)
\psline[linewidth=0.04cm](0.9409375,-1.7954688)(2.3809376,-1.7954688)
\psline[linewidth=0.04cm](2.3809376,-1.7954688)(2.3809376,-0.59546876)
\psline[linewidth=0.04cm](0.9409375,0.6045312)(2.3809376,-0.59546876)
\usefont{T1}{ptm}{m}{n}
\rput(2.87,1.9145312){\scriptsize$001$}
\usefont{T1}{ptm}{m}{n}
\rput(2.87,0.47453126){\scriptsize$011$}
\usefont{T1}{ptm}{m}{n}
\rput(2.87,-0.48546875){\scriptsize$101$}
\usefont{T1}{ptm}{m}{n}
\rput(2.87,-1.6854688){\scriptsize$111$}
\usefont{T1}{ptm}{m}{n}
\rput(0.41234374,1.9145312){\scriptsize$000$}
\usefont{T1}{ptm}{m}{n}
\rput(0.41234374,0.47453126){\scriptsize$100$}
\usefont{T1}{ptm}{m}{n}
\rput(0.41234374,-0.48546875){\scriptsize$010$}
\usefont{T1}{ptm}{m}{n}
\rput(0.41234374,-1.6854688){\scriptsize$110$}
\psline[linewidth=0.04cm,linestyle=dashed,dash=0.16cm
0.16cm](8.620937,1.8045312)(7.1809373,1.8045312)
\psline[linewidth=0.04cm,linestyle=dashed,dash=0.16cm
0.16cm](7.1809373,0.6045312)(8.620937,0.6045312)
\psline[linewidth=0.04cm,linestyle=dashed,dash=0.16cm
0.16cm](7.1809373,-0.59546876)(8.620937,-0.59546876)
\psline[linewidth=0.04cm,linestyle=dashed,dash=0.16cm
0.16cm](7.1809373,-1.7954688)(8.620937,-1.7954688)
\usefont{T1}{ptm}{m}{n}
\rput(1.6923437,-2.7){\scriptsize$EH(1,1)$}
\usefont{T1}{ptm}{m}{n}
\rput(7.8,-2.7){\scriptsize$EH(1,2)$}
\end{pspicture}
\end{center}
\caption{\label{f1}\footnotesize {Two exchanged hypercubes $EH(1,1)$ and $EH(1,2)$}}
\end{figure}
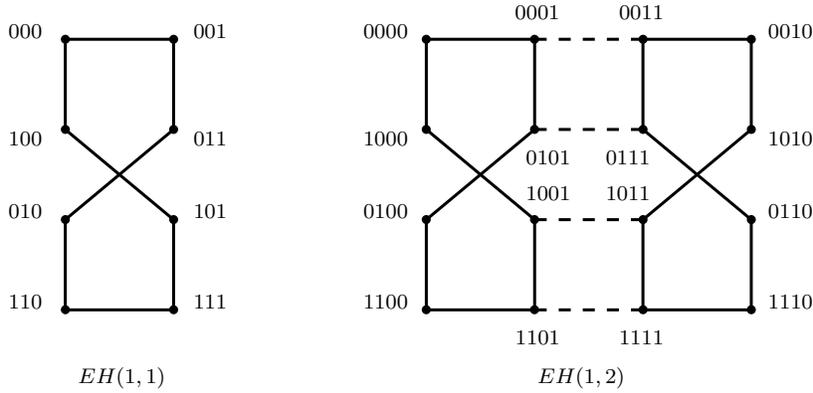

The exchanged hypercubes $EH(1, 1)$ and $EH(1, 2)$ are shown in
Figure~\ref{f1}.

From Definition~\ref{def2.2}, it is easy to see that $EH(s,t)$ can
be obtained from a hypercube $Q_{s+t+1}$ with vertex-set $V(s,t)$ by
removing all $r$-dimensional edges that link two vertices with the
last bit $0$ if $r\in I_t$ and two vertices with the last bit $1$ if
$r\in I_{s+t}-I_t$.
Thus, $EH(s, t)$ is a bipartite graph with minimum degree
$\min\{s,t\}+1$ and maximum degree $\max\{s,t\}+1$. The following
three lemmas obtained by Loh {\it et al.}~\cite{lhp05} and
Ma~\cite{m10} are very useful for our proofs.

\begin{lem}\label{lem2.3} {\rm(Loh {\it et al.}~\cite{lhp05})}\
$EH(s,t)$ is isomorphic to $EH(t,s)$.
\end{lem}

By Lemma~\ref{lem2.3}, without loss of generality, we can assume
$s\leqslant t$ in the following discussion, and so $EH(s,t)$ has the
minimum degree $s+1$. For fixed $r \in I_{s+t}$ and $i \in \{0,1\}$,
let $H^{r}_i$ denote a subgraph of $EH(s,t)$ induced by all vertices
whose the $r$-th bits are $i$.

\begin{lem}\label{lem2.4} {\rm(Loh {\it et al.}~\cite{lhp05})}\
For a fixed $r \in I_{s+t}$, $EH(s,t)$ can be decomposed into $2$
isomorphic subgraphs $H^r_0$ and $H^r_1$, which are isomorphic to
$EH(s,t-1)$ if $r\in I_{t}$ and $t\geqslant 2$, and isomorphic to
$EH(s-1,t)$ if $r\in I_{s+t}-I_{t}$ and $s\geqslant 2$. Moreover,
there are $2^{s+t-1}$ independent edges between $H^r_0$ and $H^r_1$.
\end{lem}

\begin{lem}\label{lem2.5} {\rm (Ma~\cite{m10})}\
$\kappa( EH(s, t) )= \lambda( EH(s, t) )=s+1$ for any $s$ and $t$
with $1\leqslant s\leqslant t$.
\end{lem}

\section{Main results}

In this section, we present our main results, that is, we determine
the $h$-connectivity and $h$-edge-connectivity of the exchanged
hypercube $EH(s,t)$.

\begin{lem}\label{lem3.1}
$\kappa^{(h)}(EH(s,t))\leqslant 2^h(s+1-h) $ and
$\lambda^{(h)}(EH(s,t))\leqslant 2^h(s+1-h)$ for $ h \leqslant s$.
\end{lem}

\begin{pf}
Let $X$ be a subset of vertices in $EH(s,t)$ whose the rightmost
$s+t+1-h$ bits are zeros and the leftmost $h$ bits do not care,
denoted by
 $$
 X=\{*^h0^{s+t+1-h}|\ *\in\{0,1\}\}.
 $$
Then the subgraph of $EH(s,t)$ induced by $X$ is a hypercube $Q_h$.
Let $S$ be the neighbor-set of $X$ in $EH(s,t)-X$ and $F$ the
edge-sets between $X$ and $S$. By Definition~\ref{def2.2}, $S$ has
the form
 $$
 S=\{*^h\underbrace{0^p10^{s-h-p-1}}_{s-h}0^{t+1}|\ 0\leqslant p\leqslant s-h-1,\ *\in\{0,1\}\}\cup\{*^h0^{s+t-h}1\}.
 $$
On the one hand, since every vertex of $X$ has degree $s+1$ in
$EH(s,t)$ and $h$ neighbors in $X$, it has exactly $s-h+1$ neighbors
in $S$. On the other hand, every vertex of $S$ has exactly one
neighbor in $X$. It follows that
 $$
  |S|=|F|=2^h(s+1-h).
 $$

We show that $S$ is an $h$-vertex-cut of $EH(s,t)$. Clearly, $S$ is
a vertex-cut of $EH(s,t)$ since $|X\cup S|=2^h(s+2-h)<2^{s+t+1}$.
Let $Y=EH(s,t)-(X\cup S)$ and $v$ be any vertex in $Y$. We only need
to show that the vertex $v$ has degree at least $h$ in $Y$.
In fact, it is easy to see from the formal definition of $S$ that if
$v$ is adjacent to some vertex in $S$ then it has only the form
 $$
 v=*^h\underbrace{0^p10^{s-h-p-1}}_{s-h}0^{t}1\ \ {\rm or}
 *^h0^{s-h}\underbrace{0^r10^{t-r-1}}_{t}1\  {\rm or} *^h\underbrace{0^p10^q10^{s-h-p-q-2}}_{s-h}0^{t+1} 
 $$
If $v$ has the former two forms, then $v$ has one neighbor in $S$, thus $v$ has at least
$(s+1-1=s\geqslant)\, h$ neighbors in $Y$.
If $v$ has the last form, then $s-h\geqslant 2 $ and $v$ has two neighbors in $S$. Thus, $v$ has
at least $(s+1-2=s-1>)\, h$ neighbors in $Y$.

By the arbitrariness of $v\in Y$, $S$ is an $h$-vertex-cut of
$EH(s,t)$, and so
 $$
 \kappa ^{(h)}(EH(s,t))\leqslant |S|=2^h(s+1-h)
 $$
as required.

We now show that $F$ is an $h$-edge-cut of $EH(s,t)$. Since every
vertex $v$ in $EH(s,t)-X$ has at most one neighbor in $X$, then  $v$
has at least $(s+1-1=s\geqslant)\, h$ neighbors in $EH(s,t)-X$. By
the arbitrariness of $v\in EH(s,t)-X$, $F$ is an $h$-edge-cut of
$EH(s,t)$, and so
 $$
 \lambda^{(h)}(EH(s,t))\leqslant |F|=2^h(s+1-h)
 $$
The lemma follows.
\end{pf}

\begin{cor}\label{cor3.2}
$\kappa^{(1)}(EH(1,t))=\lambda^{(1)}(EH(1,t))=2$ for $t\geqslant 1$.
\end{cor}

\begin{pf}\label{cor3.1}
On the one hand, $\kappa ^{(h)} (EH(1,t))\leqslant 2$ and
$\lambda^{(h)} (EH(1,t))\leqslant 2$ by Lemma~\ref{lem3.1} when
$s=1$. On the other hand, by Lemma~\ref{lem2.5},
$\kappa(EH(1,t))=\lambda(EH(1,t))=2$, thus $\kappa ^{(h)}(EH(1,t))
\geqslant \kappa(EH(1,t))=2$ and $\lambda ^{(h)}(EH(1,t)) \geqslant
\lambda(EH(1,t))=2$. The results hold.
\end{pf}
\vskip6pt

\begin{thm}\label{thm3.4}
For $1\leqslant s\leqslant t$ and any $h$ with $0\leqslant h\leqslant s$,
 $$
 \kappa^{(h)}(EH(s,t))=\lambda^{(h)}(EH(s,t))= 2^h(s+1-h).
 $$
\end{thm}

\begin{pf}
By Lemma~\ref{lem3.1}, we only need to prove that,
 $$
 \kappa^{(h)}(EH(s,t))= \lambda^{(h)}(EH(s,t))\geqslant 2^h(s+1-h).
 $$
We proceed by induction on $h\geqslant 0$. The theorem holds for
$h=0$ by Lemma~\ref{lem2.5}. Assume the induction hypothesis for
$h-1$ with $h\geqslant 1$, that is,
  \begin{eqnarray}\label{e3.1}
 \kappa ^{(h-1)}(EH(s,t))=\lambda ^{(h-1)}(EH(s,t))\geqslant 2^{h-1}(s+2-h).
 \end{eqnarray}

Note $h=1$ if $s=1$. By Corollary~\ref{cor3.2},
$\kappa^{(1)}(EH(1,t))=\lambda^{(1)}(EH(1,t))=2$ for any $t\geqslant
1$, the theorem is true for $s=1$. Thus, we assume $s\geqslant 2$
below.

Let $S$ be a minimum $h$-vertex-cut (or $h$-edge-cut) of $EH(s,t)$
and $X$ be the vertex-set of a minimum connected component of
$EH(s,t)-S$. Then
 $$
 |S|=\left\{
 \begin{array}{rl}
 \kappa^{(h)}(EH(s,t)) &\ {\it if}\  S\ {\it is\ a\ vertex-cut};\\
 \lambda^{(h)}(EH(s,t)) &\ {\it if}\  S\ {\it is\ an\ edge-cut}.
 \end{array}\right.
 $$
Thus, we only need to prove that
 \begin{equation}\label{e3.2}
 |S|\geqslant 2^h(s+1-h).
 \end{equation}
To the end, let $Y$ be the set of vertices in $EH(s,t)-S$ not in
$X$, and for a fixed $r\in I_{s+t}$ and each $i=0,1$, let
$$
 \begin{array}{l}
 X_i=X\cap H^r_i,\\
 Y_i=Y\cap H^r_i\ {\rm and}\\
 S_i=S\cap H^r_i,
 \end{array}
$$

Let $J=\{i\in \{0,1\}|\ X_i\ne\emptyset\}$ and $J'=\{i\in  J|\
Y_i\not=\emptyset\}$. Clearly, $0\leqslant |J'|\leqslant
|J|\leqslant 2$ and $|J'|=0$ only when $|J|=1$. We choose $r\in
I_{s+t}$ such that $|J|$ is as large as possible. For each $i\in
\{0,1\}$, we write $H_i$ for $H^r_i$ for short. We first prove the
following inequality.

\begin{equation}\label{e3.3}
|S_i|\geqslant 2^{h-1}(s+1-h)\ {\rm if}\ X_i \neq \emptyset\ {\rm
and}\ Y_i\neq \emptyset\ {\rm for}\ i\in \{0,1\}.
\end{equation}

In fact, for some $i\in\{0,1\}$, if $X_i \neq \emptyset$ and
$Y_i\neq \emptyset$, then $S_i$ is a vertex-cut (or an edge-cut) of
$H_i$. Let $u$ be any vertex in $X_i\cup Y_i$. Since $S$ is an
$h$-vertex-cut (or $h$-edge-cut) of $EH(s,t)$, $u$ has degree at
least $h$ in $EH(s,t)-S$. By Lemma~\ref{lem2.4}, $u$ has at most one
neighbor in $H_j$, where $j\ne i$. Thus, $u$ has degree at least
$h-1$ in $H_i$, which implies that $S_i$ is an $(h-1)$-vertex-cut
(or edge-cut) of $H_i$, that is,
 \begin{equation}\label{e3.4}
 |S_i|\geqslant \kappa^{(h-1)}(H_i)\ \ {\rm (or}\ |S_i|\geqslant \lambda^{(h-1)}(H_i)).
 \end{equation}

If $r \in I_{s+t}-I_{t}$, then $H_i\cong EH(s-1,t)$ by
Lemma~\ref{lem2.4}. By the induction hypothesis (\ref{e3.1}),
 $
 \kappa^{(h-1)}(H_i)=\lambda^{(h-1)}(H_i) \geqslant
2^{h-1}(s+1-h),
 $
from which and (\ref{e3.4}), we have that
 $
|S_i|\geqslant 2^{h-1}(s+1-h).
 $

If $r \in I_{t} $, then $H_i\cong EH(s,t-1)$ by Lemma~\ref{lem2.4}.

If $t \geqslant s+1 $, by the induction hypothesis (\ref{e3.1}),
 $$
\kappa^{(h-1)}(H_i)=\lambda^{(h-1)}(H_i)\geqslant 2^{h-1}(s+2-h)>
2^{h-1}(s+1-h),
$$
from which and (\ref{e3.4}), we have that $|S_i|> 2^{h-1}(s+1-h)$.

If $t=s$, then $EH(s,t-1)\cong EH(s-1,t)$ by Lemma~\ref{lem2.3}. By
the induction hypothesis (\ref{e3.1}),
 $$
 \kappa^{(h-1)}(H_i)=\lambda^{(h-1)}(H_i)\geqslant 2^{h-1}(s+1-h),
 $$
from which and (\ref{e3.4}), we have that $|S_i|\geqslant
2^{h-1}(s+1-h)$. The inequality (\ref{e3.3}) follows.

  \vskip6pt
We now prove the inequality in (\ref{e3.2}).

If $|J|=1$ then, by the choice of $J$, no matter what $r\in I_{s+t}$
is chosen, the $r$-th bits of all vertices in $X$ are the same. In
other words, the $r$-th bits of all vertices in $X$ are the same for
any $r\in I_{s+t}$, and possible different in the last bit. Thus
$|X|\leqslant 2$ and $h\leqslant 1$. By the hypothesis of
$h\geqslant 1$, we have $h=1$ and $|X|=2$. The subgraph of $EH(s,t)$
induced by $X$ is an edge in $E_1$, thus
 $$
 |S|=s+t\geqslant 2s =2^{h}(s+1-h),
 $$
as required. Assume $|J|=2$ below, that is, $X_i \neq \emptyset$ for
each $i=0,1$. In this case, $|J'|\geqslant 1$.

If  $|J'|=2$ then, for each $i=0,1$, since $X_i \neq \emptyset$ and
$Y_i\neq \emptyset$, we have that $|S_i|\geqslant 2^{h-1}(s+1-h)$ by
(\ref{e3.3}). Note that $|S|=|S_0|+|S_1|$ if $S$ is an
$h$-vertex-cut and $|S|\geqslant |S_0|+|S_1|$ if $S$ is an
$h$-edge-cut. It follows that
$$
\begin{array}{rl}
|S|&\geqslant |S_0|+|S_1| \\
&\geqslant 2\times 2^{h-1}(s+1-h)\\&=2^{h}(s+1-h),
 \end{array}
$$
as required.

If $|J'|=1$, then one of $Y_0$ and $Y_1$ must be empty. Without loss
of generality, assume $Y_1=\emptyset$ and $ Y_0 \neq \emptyset$.

Clearly, $S$ is not an $h$-edge-cut, otherwise, $|Y|<|H_0|<|X|$, a
contradiction with the minimality of $X$. Thus, $S$ is an
$h$-vertex-cut. By (\ref{e3.3}), $|S_0|\geqslant 2^{h-1}(s+1-h)$.
Since $Y_1=\emptyset$, we have
 \begin{equation}\label{e3.5}
 |X_1|=|H_1|-|S_1| \ \ {\rm and}\ \ |Y|=|H_0|-|X_0|-|S_0|.
\end{equation}
If $|S_1|< |S_0|$ then, by (\ref{e3.5}), we obtain that
$|Y|<|X_1|<|X|$, which contradicts to the minimality of $X$. Thus,
$|S_1|\geqslant |S_0|$, from which and (\ref{e3.3}) we have that
 $$
 \begin{array}{rl}
 |S|&=|S_0|+|S_1|\geqslant 2|S_0|\\
 &\geqslant 2\times2^{h-1}(s+1-h)\\
 &=2^{h}(s+1-h),
  \end{array}
$$
as required. Thus, the inequality in (\ref{e3.2}) holds, and so the
theorem follows.
\end{pf}

\begin{cor} {\rm (Ma and Zhu~\cite{m11})}\
If $1\leqslant s\leqslant t$, then
$\kappa^{(1)}(EH(s,t))=\lambda^{(1)}(EH(s,t))=2s$.
\end{cor}

A dual-cube $DC(n)$, proposed by Li and Peng~\cite{lp00} constructed
from hypercubes, preserves the main desired properties of the
hypercube. Very recently, Yang and Zhou~\cite{yz12} have determined
that $\kappa^{(h)}(DC(n))=2^n(n+1-h)$ for each $h=0,1,2$. Since
$EH(n,n)$ is isomorphic to $DC(n)$, the following result is obtained
immediately.

\begin{cor}
For dual-cube $DC(n)$,
$\kappa^{(h)}(DC(n))=\lambda^{(h)}(DC(n))=2^n(n+1-h)$ for any $h$
with $0\leqslant h\leqslant n$.
\end{cor}

\section{Conclusions}

In this paper, we consider the generalized measures of of fault
tolerance for a network, called the $h$-connectivity $\kappa^h$ and
the $h$-edge-connectivity $\lambda^h$. For the exchanged hypercube
$EH(s,t)$, which has about half edges of the hypercube $Q_{s+t+1}$,
we prove that $\kappa^{(h)}=\lambda^{(h)}= 2^h(s+1-h)$ for any $h$
with $0\leqslant h\leqslant s$ and $s\leqslant t$. The results show
that at least $2^h(s+1-h)$ vertices (resp. $2^h(s+1-h)$ edges) of
$EH(s,t)$ have to be removed to get a disconnected graph that
contains no vertices of degree less than $h$. Thus, when the
exchanged hypercube is used to model the topological structure of a
large-scale parallel processing system, these results can provide
more accurate measurements for fault tolerance of the system.

Otherwise, Ma and Liu~\cite{ml09} investigated bipancyclicity of
$EH(s,t)$. However, there are many interesting combinatorial and
topological problems, e.g., wide-diameter, fault-diameter,
panconnectivity, spanning-connectivity, which are still open for the
exchanged hypercube network.


\begin{thebibliography}{10}


\bibitem{E89}
A. H. Esfahanian, Generalized measures of fault tolerance with
application to $n$-cube networks. IEEE Transactions on Computers, 38
(11) (1989), 1586-1591.

\bibitem{EH88}
A. H. Esfahanian, S.L. Hakimi, On computing a conditional edge
connectivity of a graph. Information Processing Letters, 27
(1988),195-199.

\bibitem{LHP94}
S. Latifi, M. Hegde, M. Naraghi-Pour, Conditional connectivity
measures for large multiprocessor systems. IEEE Transactions on
Computers, 43 (1994) 218-222.

\bibitem{lx12a}
X.-J. Li and J.-M. Xu, Generalized measures of fault tolerance in
$(n,k)$-star graphs.  http://arxiv.org/abs/1204.1440, 2012.

\bibitem{lx12}
X.-J. Li and J.-M. Xu, Generalized measures of edge fault tolerance
in $(n,k)$-star graphs. Mathematical Science Letters, 1 (2) (2012),
133-138.

\bibitem{lp00}
Y. Li, S. Peng, Dual-cubes: a new interconnection network for
high-performance computer clusters. In: Proceedings of the 2000
International Computer Architecture, (2000), pp. 51-57.


\bibitem{lhp05}
P. K. K. Loh, W. J. Hsu, Y. Pan, The exchanged hypercube. IEEE
Transactions on Parallel and Distributed Systems, 16 (9) (2005),
866-874.

\bibitem{m10}
M. Ma, The connectivity of exchanged hypercubes. Discrete
Mathematics Algorithms and Applications, 2 (2) (2010), 213-220.

\bibitem{ml09}
M. Ma, B. Liu, Cycles embedding in exchanged hypercubes. Information
Processing Letters, 110 (2) (2009), 71-76.


\bibitem{m11}
M. Ma and L. Zhu, The super connectivity of exchanged hypercubes.
Information Processing letters, 111 (2011), 360-364.

\bibitem{oc93}
A. D. Oh, H. Choi, Generalized measures of fault tolerance in
$n$-cube networks. IEEE Transactions on Parallel and Distributed
Systems, 4 (1993), 702-703.

\bibitem{wz09}
M. Wan, Z. Zhang, A kind of conditional vertex connectivity of star
graphs. Applied Mathematics Letters,  22 (2009), 264-267.

\bibitem{wg98}
J. Wu, G. Guo, Fault tolerance measures for $m$-ary $n$-dimensional
hypercubes based on forbidden faulty sets. IEEE Transactions on
Computers, 47 (1998), 888-893.

\bibitem{x00c}
J.-M. Xu, On conditional edge-connectivity of graphs. Acta
Mathematae Applicatae Sinica, 16 (4) (2000), 414-419.

\bibitem{xu01}
J.-M. Xu, Topological Structure and Analysis of Interconnection
Networks. Kluwer Academic Publishers, Dordrecht/Boston/London, 2001.

\bibitem{xxz05}
J.-M. Xu, M. Xu, Q. Zhu, The super connectivity of shuffle-cubes.
Information Processing Letters, 96 (2005), 123-127.

\bibitem{ylg10}
W.-H. Yang, H.-Z. Li, X.-F, Guo, A kind of conditional fault
tolerance of $(n,k)$-star graghs.  Information Processing Letters,
110 (2010), 1007-1011.

\bibitem{yz12}
X. Yang, S. Zhou, On conditional fault tolerant of dual-cubes.
International Journal of Parallel, Emergent and Distributed Systems.
DOI:10.1080/17445760.2012.704631, 2012

\bibitem{zxh07}
Q. Zhu, J.-M. Xu, X.-M. Hou, X. Xu, On reliability of the folded
hypercubes, Information Sciences, 177 (8) (2007), 1782-1788.

\bibitem{zxl06}
Q. Zhu, J.-M. Xu, M. L\"u, Edge fault tolerance analysis of a class
of interconnection networks. Applied Mathematics and Computation,
172 (1) (2006), 111-121.

\end{thebibliography}
\end{document}